\documentstyle[12pt,amssymb,amsbsy,righttag]{amsart}
\renewcommand{\a}{\alpha}

\newcommand{\g}{\gamma}

\newcommand{\z}{\zeta}

\renewcommand{\l}{\lambda}

\renewcommand{\t}{\tau}

\newcommand{\f}{\varphi}
\renewcommand{\o}{\omega}

\renewcommand{\L}{\Lambda}

\renewcommand{\O}{\Omega}
\newcommand{\U}{\Upsilon}

\newcommand{\cR}{{\cal R}}

\newcommand{\C}{{\Bbb C}}
\newcommand{\T}{{\Bbb T}}
\newcommand{\pp}{{\Bbb P}}

\newcommand{\R}{{\Bbb R}}
\newcommand{\Z}{{\Bbb Z}}
\newcommand{\mm}{{\Bbb M}}
\newcommand{\0}{{\Bbb O}}

\newcommand{\m}{{\boldsymbol m}}
\newcommand{\bS}{{\boldsymbol S}}

\newcommand{\df}{\stackrel{\mathrm{def}}{=}}
\newcommand{\dist}{\operatorname{dist}}
\newcommand{\Ker}{\operatorname{Ker}}

\newcommand{\rank}{\operatorname{rank}}
\newcommand{\const}{\operatorname{const}}

\newcommand{\eeq}{\end{equation}}
\newcommand{\beq}{\begin{equation}}
\newcommand{\bay}{\begin{eqnarray}}
\newcommand{\ey}{\end{eqnarray}}

\newcommand{\be}{\infty}

\newcommand{\bl}{\blacksquare}
\newcommand{\ess}{\operatorname{ess}}
\newcommand{\ind}{\operatorname{ind}}

\newcommand{\Range}{\operatorname{Range}}

\newcommand{\Pf}{{\bf Proof. }}

\newcommand{\ov}{\overline}

\newtheorem{thm}{\hspace{\parindent}Theorem}[section]

\newtheorem{lem}[thm]{\hspace{\parindent}Lemma}

\begin{document}
\newcommand{\bs}{\boldsymbol}
\newcommand{\diag}{\operatorname{diag}}
\renewcommand{\theequation}{\thesection.\arabic{equation}}

\author{R.B. Alexeev and V.V. Peller}
\thanks{The second author is partially supported by NSF grant DMS 9970561}

\title{Unitary interpolants and factorization \\indices of matrix functions}

\maketitle

\begin{abstract}
For an $n\times n$ bounded matrix function $\Phi$ we study unitary 
interpolants $U$, i.e., unitary-valued functions $U$
such that $\hat U(j)=\hat\Phi(j)$, $j<0$. We are looking for unitary interpolants
$U$ for which the Toeplitz operator $T_U$ is Fredholm. 
We give a new approach based on superoptimal singular values and
thematic factorizations. We describe Wiener--Hopf factorization
indices of $U$ in terms of superoptimal singular values of $\Phi$ and thematic 
indices of $\Phi-F$, where $F$ is a superoptimal approximation of $\Phi$
by bounded analytic matrix functions. The approach essentially relies on
the notion of a monotone thematic factorization introduced in [AP]. In the last 
section we discuss hereditary properties of unitary interpolants. In particular, for 
matrix functions $\Phi$ of class $H^\be+C$ we study unitary interpolants $U$ of class 
$QC$.
\end{abstract}

\setcounter{equation}{0}
\section{\bf Introduction}

\

A Hankel operator defined on the Hardy class $H^2(\C^n)$ of $\C^n$-valued functions
has infinitely many different symbols. If $\Phi$ is a symbol of such a Hankel
operator, then $\Phi-Q$ is a symbol of the same Hankel operator for any
bounded analytic matrix function $Q$. A natural problem arises for a Hankel
operator to choose a symbol that satisfies certain nice properties. For example,
an important problem is to choose a symbol that has minimal $L^\be$-norm.
In this paper we consider another important problem to choose a symbol that takes 
unitary values (such symbols are called {\it unitary-valued}).

Recall that  for a bounded $\mm_{m,n}$-valued function $\Phi$ 
(we denote by $\mm_{m,n}$ the space of $m\times n$ matrices)
on the unit circle $\T$ the {\it Hankel operator} $H_\Phi:H^2(\C^n)\to H^2_-(\C^m)$
with symbol $\Phi$ is defined by
$$
H_\Phi f\df\pp_-\Phi f,\quad f\in H^2(\C^n),
$$
where $\pp_-$ is the orthogonal projection onto 
$H^2_-(\C^m)\df L^2(\C^m)\ominus H^2(\C^m)$. 

Certainly, when we discuss the problem of finding unitary-valued symbols we
have to assume that $m=n$, i.e., $\Phi\in L^\be(\mm_{n,n})$. If $U$ is a 
symbol of $H_\Phi$, then the Fourier coefficients of $U$ satisfy
$$
\hat U(j)=\hat\Phi(j),\quad j<0.
$$ 
A unitary-valued matrix function $U$ satisfying this condition 
is called a {\it unitary interpolant} of $\Phi$. 

A matrix analog of Nehari's theorem says that
\beq
\label{1.1}
\|H_\Phi\|=\dist_{L^\be}\big(\Phi,H^\be(\mm_{m,n})\big)
\end{equation}
(see [Pa]). Here for a bounded $\mm_{m,n}$-valued function $F$ we use the notation
$$
\|F\|_{L^\be}\df\ess\sup_{\z\in\T}\|F(\z)\|_{\mm_{m,n}},
$$
where the norm of a matrix in $\mm_{m,n}$ is its operator norm from $\C^n$ to $\C^m$.

Recall that by Hartman's theorem, $H_\Phi$ is compact if and only if \linebreak
$\Phi\in(H^\infty+C)(\mm_{m,n})$ (see e.g., [N]), where
$$
H^\be+C\df\{f+g:~f\in H^\be,~g\in C(\T)\},
$$
while the essential norm $\|H_\Phi\|_{\text e}$ (i.e., the distance from $H_\Phi$ 
to the set of compact operators) can be computed as follows: 
$$
\|H_\Phi\|_{\text e}=\dist_{L^\be}\big(\Phi,(H^\be+C)(\mm_{m,n})\big)
$$
(see e.g., [Sa] for the proof of this formula for scalar functions, in the
matrix-valued case the proof is the same).

It follows from (\ref{1.1}) that for a matrix function $\Phi\in L^\be(\mm_{n,n})$
to have a unitary interpolant it is necessary that $\|H_\Phi\|\le 1$.

If $\f$ is a scalar function such that $\|H_\f\|\le1$, then a unitary interpolant
of $\f$ exists if $H_\f$ has two different symbols whose $L^\be$-norms are at most 1
(the Adamyan--Arov--Krein theorem [AAK2]). 

Dym and Gohberg studied in [DG1] the problem of finding unitary interpolants for
matrix functions with entries in a Banach algebra $X$ of continuous functions
on $\T$ that satisfy certain axioms (their axioms are similar to axioms
(A1)--(A4) in \S 5). They showed that if $\Phi$ belongs to the space $X(\mm_{n,n})$ 
of $\mm_{n,n}$-valued functions with entries in $X$ and $\|H_\Phi\|\le1$,
then $\Phi$ has a unitary interpolant $U$ of the same class $X(\mm_{n,n})$. Moreover, 
the negative indices of a Wiener--Hopf factorization of $U$
are uniquely determined by $\Phi$ while the nonnegative
indices can be arbitrary. Note, however, that we were not able to understand
their argument in the proof of Lemma 4.10. In [DG2] the problem of finding
unitary interpolants was studied in a more general situation. Another approach
to this problem was given by Ball in [B].

In this paper we propose a different approach based on so-called thematic 
factorizations introduced in [PY1]. We state our results in different terms.
We obtain information about the Wiener--Hopf indices of unitary interpolants
of $\Phi$ in terms of the so-called superoptimal singular values of $\Phi$
and indices of thematic factorizations of $\Phi-F$, where $F$ is 
a superoptimal approximation of $\Phi$ by bounded analytic matrix functions 
(see \S 2).

We are interested in unitary interpolants $U$ of $\Phi$ such that the
Toeplitz operator $T_U$ is {\it Fredholm} (i.e., invertible modulo the
compact operators).

Recall that for $\Psi\in L^\be(\mm_{m,n})$ the {\it Toeplitz operator}
$T_\Psi:H^2(\C^n)\to H^2(\C^m)$ is defined by
$$
T_\Psi f=\pp_+\Psi f,\quad f\in H^2(\C^n),
$$
where $\pp_+$ is the orthogonal projection onto $H^2(\C^m)$.

By Simonenko's theorem [Si] (see also [LS]), the symbol $\Psi\in L^\be(\mm_{n,n})$ 
of a Fredholm Toeplitz operator $T_\Psi$ admits a {\it Wiener--Hopf factorization}:
$$
\Psi=Q_2^*\left(\begin{array}{cccc}z^{d_0}&0&\cdots&0\\
0&z^{d_1}&\cdots&0\\\vdots&\vdots&\ddots&\vdots\\0&0&\cdots&z^{d_{n-1}}\end{array}
\right)Q_1,
$$
where $d_0,\cdots,d_{n-1}\in\Z$ are {\it Wiener--Hopf indices}, and $Q_1$
and $Q_2$ are functions invertible in $H^2(\mm_{n,n})$.
Here we start enumeration with 0 for technical reasons.

Clearly, we can always arrange the indices in the nondecreasing order:
$$
d_0\le d_1\le\cdots\le d_{n-1}
$$
in which case the indices $d_j$ are uniquely determined by the function $\Psi$.
Uniqueness follows easily from the following well-known identity 
\beq
\label{1.2}
\dim\Ker T_\Psi=\sum_{\{j:d_j<0\}}-d_j
\end{equation}
applied to the matrix functions $z^j\Psi$, $j\in\Z$.

Let $U$ be a unitary-valued function on $\T$. It is well known that the Toeplitz 
operator $T_U$ is Fredholm if and only if 
\beq
\label{1.3}
\|H_U\|_{\text e}<1,\quad\mbox{and}\quad\|H_{U^*}\|_{\text e}<1.
\end{equation}
Indeed, suppose that $T_U$ is Fredholm. Clearly, 
$\|H_{z^jU}\|_{\text e}=\|H_U\|_{\text e}$ for any $j\in\Z$. Multiplying $U$ by
$z^N$ if necessary, we may assume that $\Ker T_U$ is trivial, and so $T_U$
is left-invertible. It follows now from the obvious equality
$$
\|H_Uf\|_2^2+\|T_Uf\|_2^2=\|f\|^2_2,\quad f\in H^2(\C^n),
$$
that $\|H_U\|<1$, and so $\|H_U\|_{\text e}<1$. Applying the same reasoning
to $U^*$, we find that $\|H_{U^*}\|_{\text e}<1$.

Suppose now that (\ref{1.3}) holds. It is easy to see that
$$
T^*_UT_U=I-H^*_UH_U,\quad T_UT^*_U=I-H^*_{U^*}H_{U^*},
$$
Since $\|H_U\|_{\text e}<1$ and $\|H_{U^*}\|_{\text e}<1$, it follows that
$T^*_UT_U$ and $T_UT^*_U$ are invertible modulo the compact operators which
implies that $T_U$ is Fredholm.

Obviously, if $U$ is a unitary interpolant of $\Phi$, then $H_U=H_\Phi$. Therefore
for $\Phi$ to have a unitary interpolant $U$ with Fredholm $T_U$ it is necessary
that $\|H_\Phi\|_{\text e}<1$. Throughout this paper we assume that this condition
is satisfied.

We conclude the introduction with the following fact.

\begin{lem}
\label{t1.1}
Let $U$ be a unitary-valued matrix function in $L^\be(\mm_{n,n})$
such that {\em$\|H_U\|_{\text e}<1$}. 
Suppose that $U$ admits a Wiener--Hopf factorization of the form
$$
U=Q_2^*\L Q_1,
$$
where $Q_1$ and $Q_2$ are invertible in $H^2(\mm_{n,n})$ and $\L$
is a diagonal matrix function with diagonal terms $z^{d_j}$.
Then $T_U$ is Fredholm and {\em
$$
\|H_U\|_{\text e}=\|H_{U^*}\|_{\text e}.
$$}
\end{lem} 

\Pf Multiplying $U$ by $\bar z^N$ for a sufficiently large $N$ 
if necessary, we may assume that all exponents $d_j$ are nonpositive. 
Let us show that in this case $T_U$ has dense range. Suppose that $g\in H^2(\C^n)$ and 
$g\perp\Range T_U$. Let $f$ be a polynomial in $H^2(\C^n)$. We have
$$
0=(T_Uf,g)=(Uf,g)=(Q_1f,\L^*Q_2g).
$$
Clearly, $\L^*\in H^\be(\mm_{n,n})$ and the set of functions of the form $Q_1f$ is 
dense in $H^2(\C^n)$. Hence, $\L^*Q_2g=0$, and so $g=0$ which proves that $T_U$
has dense range.

In [Pe3] it was proved that if $T_U$ has dense range in $H^2(\C^n)$,
then the operator $H^*_{U^*}H_{U^*}$ is unitarily equivalent to the restriction of
the operator $H^*_UH_U$ to the subspace
$$
H^2(\C^n)\ominus\{f\in H^2(\C^n):~\|H_Uf\|_2=\|f\|_2\}
$$
(see also [PK] where this was proved in the scalar case). Clearly, the condition
$\|H_U\|_{\text e}<1$ implies that the above subspace has finite 
codimension. Therefore $\|H_{U^*}\|_{\text e}=\|H_U\|_{\text e}<1$ and as we
have already observed, this implies that $T_U$ is Fredholm. $\bl$

To simplify the notation, for a matrix function $\Phi$ and a space $X$ of functions on 
$\T$ we can write $\Phi\in X$ when all entries of $\Phi$ belong to $X$, if this does
not lead to confusion.

In \S 2 we collect necessary information on superoptimal approximation and
thematic (as well as partial thematic) factorizations. In particular, we define
the important notion of a monotone (partial) thematic factorization that was 
introduced in [AP] and state the theorem on the invariance of indices of
monotone (partial) thematic factorizations obtained in [AP].
  
In \S 3 we study unitary interpolants of matrix functions $\Phi\in L^\be(\mm_{n,n})$
satisfying $\|H_\Phi\|_{\text e}<1$. We state the main results in terms of the 
superoptimal singular values of $\Phi$ and the indices of a monotone (partial) 
thematic factorization of $\Phi-F$, where $F$ is a best approximation of $\Phi$
by bounded analytic functions.

Finally, we show in \S 4 how to apply the results of
\S 3 to study unitary interpolants of class $X$ for various function
spaces $X$. In particular we show in \S 4 that if $\Phi\in(H^\be+C)(\mm_{n,n})$,
then all unitary interpolants $U$ of $\Phi$ satisfying the condition
$\|H_{U^*}\|_{\text e}<1$ belong to the class
$$
QC\df\{f\in H^\be+C:~\bar f\in H^\be+C\}.
$$

\

\setcounter{equation}{0}
\section{\bf Superoptimal approximation and thematic indices}

\

This section is an introduction to superoptimal approximation and thematic
factorizations. We refer the reader to [PY1], [PY2], [PT], and [AP]  for more detailed 
information. 

It is a well-known fact [Kh] that if $\f$ is a scalar function of class $H^\be+C$, then
there exists a unique function $f\in H^\be$ such that
$$
\dist_{L^\be}(\f,H^\be)=\|\f-f\|_{L^\be}.
$$
For a matrix function $\Phi\in L^\be(\mm_{m,n})$ we say that a matrix function 
$F\in H^\be(\mm_{m,n})$ is a {\it best approximation of $\Phi$ by bounded analytic 
matrix functions} if
$$
\|\Phi-F\|_{L^\be}=\dist_{L^\be}(\Phi,H^\be(\mm_{m,n}).
$$
However, consideration of diagonal matrices makes it obvious that for matrix functions 
we can have uniqueness only in exceptional cases.

To define a superoptimal approximation of a matrix function 
$\Phi\in L^\be(\mm_{m,n})$ by bounded analytic matrix functions, we
consider the following sets:
$$
\O_0=\{F\in H^\be(\mm_{m,n})
:~F~\mbox{minimizes}~ t_0\df\ess\sup_{\z\in\T}\|\Phi(\z)-F(\z)\|\};
$$
$$
\O_j=\{F\in \O_{j-1}:~F~\mbox{minimizes}~ 
t_j\df\ess\sup_{\z\in\T}s_j(\Phi(\z)-F(\z))\}.
$$
Recall that for a matrix $A\in\mm_{m,n}$ the $j$th singular value $s_j(A)$
is defined by
$$
s_j(A)=\inf\{\|A-R\|:~\rank R\le j\},\quad j\ge0.
$$ 

Functions in $F\in\O_{\min\{m,n\}-1}$ are called {\it superoptimal approximations
of $\Phi$ by analytic functions}, or superoptimal solutions of Nehari's problem.
The numbers $t_j$ are called {\it superoptimal singular values of} $\Phi$.

In [PY1] it was shown that any matrix function $\Phi\in(H^\be+C)(\mm_{m,n})$
has a unique superoptimal approximation by bounded analytic matrix functions.
Another approach to the uniqueness of a superoptimal approximation was found later 
in [T]. Those results were extended in [PT] to the case of matrix functions
$\Phi\in L^\be(\mm_{m,n})$ such that the essential norm $\|H_\Phi\|_{\text e}$
of the Hankel operator $H_\Phi$ is less than the smallest nonzero superoptimal 
singular value of $\Phi$.

A matrix function $\Phi\in L^\be(\mm_{m,n})$ is called {\it badly approximable}
if
$$
\dist_{L^\be}(\Phi,H^\be(\mm_{m,n}))=\|\Phi\|_{L^\be}.
$$
It is called {\it very badly approximable} if the zero matrix function is a 
superoptimal approximation of $\Phi$.

Recall that a nonzero scalar function $\f\in H^\be+C$ is badly approximable if and 
only if it has constant modulus almost everywhere on $\T$, belongs to $QC$, and 
$\ind T_\f>0$. For continuous $\f$ this was proved in [Po]
(see also [AAK1]). For the general case see [PK]. More generally, if $\f$ is a scalar
function in $L^\be$ such that $\|H_\f\|_{\text e}<\|H_\f\|$, then $\f$ is badly
approximable if $|\f|$ is constant almost everywhere on $\T$, $T_\f$ is Fredholm,
and $\ind T_\f>0$.

Let us now define a thematic matrix function.
Recall that a function \linebreak$F\in H^\be(\mm_{m,n})$ is called {\it inner} if
$F^*(\z)F(\z)=I_n$ almost everywhere on $\T$ ($I_n$ stands for the identity matrix
in $\mm_{n,n}$). $F$ is called {\it outer} if $FH^2(\C^n)$ is dense in $H^2(\C^m)$.
Finally, $F$ is called {\it co-outer} if the transposed function $F^{\text t}$
is outer. 

An $n\times n$ matrix function $V$, $n\ge2$, is called {\it thematic} (see [PY1])
if it is unitary-valued and has the form
$$
V=\left(\begin{array}{cc}\bs{v}&\ov{\Theta}\end{array}\right),
$$
where the matrix functions $\bs{v}\in H^\be(\C^n)$ and $\Theta\in H^\be(\mm_{n,n-1})$
are both inner and co-outer. Note that if $V$ is a thematic function, then
all minors of $V$ on the first column (i.e., all minors of $V$ of an arbitrary size 
that involve the first column of $V$) belong to $H^\be$ ([PY1]). If $n=1$ the thematic
functions are constant functions whose modulus is equal to 1.

Let now $\Phi\in L^\be(\mm_{m,n})$ be a matrix function such that
$\|H_\Phi\|_{\text e} <\|H_\Phi\|$. It follows from the results of [PT] that $\Phi$
is badly approximable if and only if it admits a representation
$$
\Phi=W^*\left(\begin{array}{cc}su&0\\0&\Psi\end{array}\right)V^*,
$$
where $s>0$, $V$ and $W^{\text t}$ are thematic functions, $u$ is a scalar unimodular
function (i.e., $|u(\z)|=1$ for almost all $\z\in\T$) such that $T_u$ is Fredholm,
$\ind T_u>0$, and $\|\Psi\|_{L^\be}\le s$. Note that in this case $s$ must be
equal to $\|H_\Phi\|$. In the case  $\Phi\in(H^\be+C)(\mm_{m,n})$ this was proved
earlier in [PY1].

Suppose now that $\|H_\Phi\|_{\text e}$ is less than the smallest nonzero superoptimal
singular value of $\Phi$ and let $m\le n$. It was proved in [PT] that 
$\Phi$ is very badly approximable if and only if
$\Phi$ admits a {\it thematic factorization}, i.e.,
\beq
\label{2.1}
\Phi = W_0^* W_1^* \cdots W^*_{m-1} D V^*_{m-1} V^*_{m-2} \cdots V^*_0
\end{equation}
where an $m \times n$ matrix function $D$ has the form
$$
D=\left( \begin{array}{ccccccc}
s_0 u_0 & 0 & \cdots & 0 & 0 & \cdots & 0 \\
0 & s_1 u_1 & \cdots & 0 & 0 & \cdots & 0\\
\vdots & \vdots & \ddots & \vdots & \vdots & \ddots & \vdots\\
0 & 0 & \cdots & s_{m-1} u_{m-1} & 0 & \cdots & 0 \end{array}
\right),
$$
$u_0, \cdots, u_{m-1}$ are unimodular scalar functions such that the operators
$T_{u_j}$ are Fredholm, $\ind T_{u_j}>0$, $\{s_j\}_{0\le j\le m-1}$ is a 
nonincreasing sequence of nonnegative numbers,
\beq
\label{2.2}
W_j=\left(\begin{array}{cc}I_j & 0 \\ 0 & \breve{W}_j \end{array} \right),\quad
V_j = \left( \begin{array}{cc} I_j & 0 \\ 0 & \breve{V}_j \end{array} \right),
\quad 1 \leq j \leq m-1,
\end{equation}
and $W^{\text t}_0, \breve{W}^{\text t}_j, V_0, \breve{V}_j$ are thematic matrix 
functions, $1 \le j\le m-1$. Moreover, in this case the $s_j$ are the superoptimal
singular values of $\Phi$: $s_j=t_j$, $0\le j\le m-1$.
The indices $k_j$ of the thematic factorization (\ref{2.1}) ({\it thematic indices}) 
are defined in case $t_j\neq0$: $k_j=\ind T_{u_j}$. Note that this result was
established earlier in [PY1] in the case $\Phi\in(H^\be+C)(\mm_{m,n})$.

It also follows from the results of [PT] that if $r\le\min\{m,n\}$ is such that
$t_r<t_{r-1}$, $\|H_\Phi\|_{\text e}<t_{r-1}$
and $F$ is a matrix function in $\O_{r-1}$, then
$\Phi-F$ admits a factorization
\beq
\label{2.3}
\Phi-F=
W_0^*\cdots W^*_{r-1} 
\left( \begin{array}{ccccccc}
t_0 u_0 & 0 & \cdots & 0 & 0 \\
0 & t_1 u_1 & \cdots & 0 & 0 \\
\vdots & \vdots & \ddots & \vdots & \vdots\\
0 & 0 & \cdots & t_{r-1}u_{r-1} & 0\\
0&0&\cdots&0&\Psi\end{array}
\right)
V^*_{r-1}\cdots V^*_0,
\end{equation}
in which the $V_j$ and $W_j$ have the form (\ref{2.2}), the
$W^{\text t}_0, \breve{W}^{\text t}_j, V_0, \breve{V}_j$ are thematic 
matrix functions, the $u_j$ are unimodular functions such that $T_{u_j}$ is
Fredholm and \linebreak$\ind T_{u_j}>0$, and $\|\Psi\|_{L^\be}\le t_r$ and
$\|H_\Psi\|<t_r$.

Factorizations of the form (\ref{2.3}) with a nonincreasing sequence
$\{t_j\}_{0\le j\le r-1}$ are called {\it partial thematic factorizations}.
If $\Phi-F$ admits a partial thematic factorization of the form (\ref{2.3}),
then $t_0,t_1,\cdots,t_{r-1}$ are the largest $r$ superoptimal singular values of
$\Phi$, and so they do not depend on the choice of a partial thematic factorization.

It was observed in [PY1] that the indices of a thematic factorization are not 
uniquely determined by the matrix function but may depend on the choice of
a thematic factorization. On the other hand it follows from the results of [PT] that 
under the conditions $\|H_\Phi\|_{\text e}<\|H_\Phi\|$ and
$t_{r-1}>\|H_\Phi\|_{\text e}$, the sum of the indices of a (partial) thematic
factorization that correspond to all superoptimal singular values equal to a
positive specific value is uniquely determined by the function $\Phi$ itself.
Note that for $H^\be+C$ matrix functions this was proved earlier in [PY2].

In [AP] the notion of a monotone (partial) thematic factorization was introduced.
It plays a crucial role in this paper. A (partial) thematic factorization is called 
{\it monotone} if for any positive superoptimal singular value $t$ the thematic
indices $k_r,k_{r+1},\cdots,k_s$ that correspond to all superoptimal singular
values equal to $t$ satisfy
$$
k_r\ge k_{r+1}\ge\cdots\ge k_s.
$$

It was shown in [AP] that under the conditions $\|H_\Phi\|_{\text e}<\|H_\Phi\|$ and
\linebreak$t_{r-1}>\|H_\Phi\|_{\text e}$, $\Phi-F$ admits a monotone partial thematic 
factorization of the form (\ref{2.3}). Moreover, it was established in [AP] that
the indices of a monotone partial thematic factorization are uniquely determined by 
$\Phi$ and do not depend on the choice of a partial thematic factorization.

\

\setcounter{equation}{0}
\section{\bf Unitary interpolants and Wiener--Hopf indices}

\

We study in this section the problem of finding unitary interpolants 
for a square matrix function. In other words, given 
$\Phi\in L^\be(\mm_{n,n})$, the problem is to find a unitary-valued function
$U$ such that $\Phi-U\in H^\be(\mm_{n,n})$. For such a problem to be solvable, 
the function $\Phi$ must satisfy the obvious necessary condition
\beq
\label{3.1}
\|H_\Phi\|=\dist_{L^\be}\big(\Phi,H^\be(\mm_{n,n})\big)\le1.
\end{equation}
However, it is well known that this condition is not sufficient even for
scalar functions (see [AAK1]).

We assume that the matrix function $\Phi$ satisfies the following condition
\beq
\label{3.2}
\|H_\Phi\|_{\text e}=\dist_{L^\be}\big(\Phi,(H^\be+C)(\mm_{n,n})\big)<1.
\end{equation}
and we are interested in unitary interpolants $U$ of $\Phi$ such that
the Toeplitz operator $T_U$ is Fredholm.
We show that under conditions (\ref{3.1}) and (\ref{3.2}) the function $\Phi$ has
such a unitary interpolant. Clearly, (\ref{3.2}) is satisfied if 
$\Phi\in(H^\be+C)(\mm_{n,n})$. 

As we have already observed in the introduction, a unitary interpolant
$U$ of a matrix function $\Phi$ satisfying (\ref{3.2}) is
the symbol of a Fredholm Toeplitz operator if and only if
\beq
\label{3.3}
\|H_{U^*}\|_{\text e}<1.
\end{equation} 
We describe all possible indices of Wiener--Hopf factorizations of such unitary 
interpolants of $\Phi$ in terms of superoptimal singular values and indices of a 
monotone partial thematic factorization of $\Phi-Q$, where $Q$ is a best approximation
of $\Phi$ by bounded analytic matrix functions.  

Let us first state the results.
As we have explained in \S 2, if $\Phi$ satisfies conditions
(\ref{3.1}) and (\ref{3.2}) and $Q\in H^\be(\mm_{n,n})$ is a best approximation
of $\Phi$ by bounded analytic matrix functions, then $\Phi-Q$ admits a monotone 
partial thematic factorization of the form
$$
\Phi-Q=W_0^*\cdots W_{r-1}^*
\left(\begin{array}{ccccc}u_0&0&\cdots&0&0\\
0&u_1&\cdots&0&0\\
\vdots&\vdots&\ddots&\vdots&\vdots\\
0&0&\cdots&u_{r-1}&0\\
0&0&\cdots&0&\Psi
\end{array}\right)V^*_{r-1}\cdots V_0^*,
$$
where $\|\Psi\|_{L^\be}\le1$ and $\|H_\Psi\|<1$ (here $r$ is the number of 
superoptimal singular values of $\Phi$ equal to 1; it may certainly happen that $r=0$
in which case $\Psi=\Phi-Q$). We denote by $k_j$, $0\le j\le r-1$,
the thematic indices of the above factorization.

\begin{thm}
\label{t3.1}
Let $\Phi$ be a matrix function in $L^\be(\mm_{n,n})$  such that
{\em$\|H_\Phi\|_{\text e}<1$}. Then $\Phi$ has a unitary interpolant $U$ satisfying
{\em$\|H_{U^*}\|_{\text e}<1$} if and only if $\|H_\Phi\|\le1$.
\end{thm}

If $U$ is a unitary interpolant of a matrix function $\Phi$ and conditions 
(\ref{3.1}), (\ref{3.2}), and (\ref{3.3}) hold, we denote by $d_j$, 
$0\le j\le n-1$, the Wiener--Hopf factorization indices of $U$
arranged in the nondecreasing order: $d_0\le d_1\le\cdots\le d_{n-1}$. 
Recall that the indices $d_j$ are uniquely determined by the function $U$.

\begin{thm}
\label{t3.2}
Let $\Phi$ be a matrix function in $L^\be(\mm_{n,n})$ such that $\|H_\Phi\|\le1$
and {\em$\|H_\Phi\|_{\text e}<1$}. Let $r$ be the number of superoptimal singular 
values of $\Phi$ equal to 1. Then each unitary interpolant $U$ of $\Phi$ satisfying 
{\em$\|H_{U^*}\|_{\text e}<1$} has precisely $r$ negative Wiener--Hopf indices. 
Moreover, $d_j=-k_j$, $0\le j\le r-1$.
\end{thm}

In particular, Theorem \ref{t3.2} says that the negative Wiener--Hopf indices of
a unitary interpolant $U$ of $\Phi$ that satisfies (\ref{3.3}) are uniquely 
determined by $\Phi$.

\begin{thm}
\label{t3.3}
Let $\Phi$ and $r$ satisfy the hypotheses of Theorem {\em\ref{t3.2}}.
Then for any sequence of integers $\{d_j\}_{r\le j\le n-1}$ satisfying
$$
0\le d_r\le d_{r+1}\le\cdots\le d_{n-1}
$$
there exists a unitary interpolant $U$ of $\Phi$ such that 
{\em$\|H_{U^*}\|_{\text e}<1$} and the nonnegative Wiener--Hopf factorization indices 
of $U$ are $d_r,d_{r+1},\cdots,d_{n-1}$.
\end{thm}

Note that Theorem \ref{t3.1} follows immediately from Theorem \ref{t3.3}.

\begin{thm}
\label{t3.4}
Let $\Phi$ satisfy the hypotheses of Theorem {\em\ref{t3.2}}.
Then $\Phi$ has a unique unitary interpolant 
if and only if all superoptimal singular values of $\Phi$ are equal to 1.
\end{thm}

{\bf Proof of Theorem \ref{t3.2}.} Suppose that $U$ is a unitary interpolant 
of $\Phi$ that satisfies (\ref{t3.3}). Put $G=U-\Phi\in H^\be(\mm_{n,n})$.

Let $\kappa\ge0$. Let us show that if $f\in H^2(\C^n)$, then
$\|H_{z^\kappa\Phi}f\|_2=\|f\|_2$ if and only if $f\in\Ker T_{z^\kappa U}$.

Indeed, if $f\in\Ker T_{z^\kappa U}$, then $z^\kappa Uf\in H^2_-(\C^n)$. 
It follows that
$$
\|H_{z^\kappa\Phi}f\|_2=\|H_{z^\kappa(\Phi+G)}f\|_2=
\|\pp_-(z^\kappa Uf)\|_2=\|z^\kappa Uf\|_2=\|f\|_2.
$$
Conversely, suppose that $\|H_{z^\kappa\Phi}f\|_2=\|f\|_2$. 
Then $\|H_{z^\kappa U}f\|_2=\|f\|_2$.
Hence, $z^\kappa Uf\in H^2_-(\C^n)$, and so  $f\in\Ker T_{z^\kappa U}$.

It follows from formula (\ref{1.2}) that
$$
\dim\Ker T_{z^\kappa U}=\sum_{\{j\in[0,n-1]:-d_j>\kappa\}}-d_j-\kappa.
$$

By Theorem 4.8 of [AP],
$$
\dim\{f\in H^2(\C^n):\|H_{z^\kappa\Phi}f\|_2=t_0\|f\|_2\}=
\sum_{\{j\in[0,r-1]:k_j>\kappa\}}k_j-\kappa.
$$
Hence,
\beq
\label{3.4}
\sum_{\{j\in[0,n-1]:-d_j>\kappa\}}-d_j-\kappa=
\sum_{\{j\in[0,r-1]:k_j>\kappa\}}k_j-\kappa,\quad \kappa\ge0.
\end{equation}

It is easy to see from (\ref{3.4}) that $U$ has $r$ negative Wiener--Hopf
factorization indices and $d_j=-k_j$ for $0\le j\le r-1$. $\bl$

To prove Theorem \ref{t3.3}, we need two lemmas.

\begin{lem}
\label{t3.5}
Let $\Psi\in L^\be(\mm_{m,m})$ and $\|H_\Psi\|<1$. Then for any nonnegative
integers $d_j$, $0\le j\le m-1$, there exists a unitary interpolant $U$ of $\Psi$ 
that admits a representation
\beq
\label{3.5}
U=W_0^*\cdots W^*_{m-1}
\left(\begin{array}{cccc}
u_0&0&\cdots&0\\
0&u_1&\cdots&0\\
\vdots&\vdots&\ddots&\vdots\\
0&0&\cdots&u_{m-1}\\
\end{array}\right)V^*_{m-1}\cdots V_0^*,
\end{equation}
where
$$
W_j=\left(\begin{array}{cc}I_j&0\\0&\breve W_j\end{array}\right),\quad
V_j=\left(\begin{array}{cc}I_j&0\\0&\breve V_j\end{array}\right),\quad
1\le j\le m-1,
$$
$V_0$, $W_0^{\text t}$, $\breve V_j$, $\breve W^{\text t}_j$ are thematic 
matrix functions, and $u_0,\cdots,u_{m-1}$ 
are unimodular functions such that the Toeplitz operators $T_{u_j}$ are Fredholm and 
$$
\ind T_{u_j}=-d_j,\quad0\le j\le m-1.
$$
\end{lem}

{\bf Proof of Lemma \ref{t3.5}.} We argue by induction on $m$. Assume first that
$m=1$. Let $\psi$ be a scalar function in $L^\be$ such that $\|H_\psi\|<1$.
Without loss of generality we may assume
that $\|\psi\|_\be<1$. Consider the function $\bar z^{d_0+1}\psi$. 
Clearly, $\|H_{ \bar z^{d_0+1}\psi}\|<1$. It is easy to see that there exists
$c\in\R$ such that
$$
\|H_{\bar z^{d_0+1}\psi+c\bar z}\|=1.
$$
Since $H_{c\bar z}$ has finite rank, it is easy to see that
$$
\|H_{\bar z^{d_0+1}\psi+c\bar z}\|_{\text e}=\|H_{\bar z^{d_0+1}\psi}\|_{\text e}
<1=\|H_{\bar z^{d_0+1}\psi+c\bar z}\|.
$$
As we have mentioned in \S 2, $\bar z^{d_0+1}\psi+c\bar z$ has a unique best
approximation $g$ by $H^\be$ functions, the error function
$u=\bar z^{d_0+1}\psi+c\bar z-g$ is unimodular, $T_u$ is Fredholm, and
$\ind T_u>0$. On the other hand
$$
\dist_{L^\be}(zu,H^\be)=\|H_{\bar z^{d_0}\psi+c-zg}\|=\|H_{\bar z^{d_0}\psi}\|<1,
$$
and so $zu$ is not badly approximable. Hence, $\ind T_{zu}\le0$. It follows
that $\ind T_u=1$. We have
$$
z^{d_0+1}u=\psi+cz^{d_0}-z^{d_0+1}g.
$$
Put $u_0=z^{d_0+1}u$. Then $\psi-u_0=z^{d_0+1}g-cz^{d_0}\in H^\be$. Hence, $u_0$ is a 
unitary interpolant of $\psi$ and $\ind T_{u_0}=-d_0-1+\ind T_u=-d_0$.

Suppose now that the lemma is proved for $(m-1)\times(m-1)$ matrix functions.
Again, without loss of generality we may assume that $\|\Psi\|_{L^\be}<1$. Then
\linebreak$\|H_{\bar z^{d_0+1}\Psi}\|<1$. 
As in the scalar case, there exists $c\in\R$ such that
$$
\|H_{\bar z^{d_0+1}\Psi+c\bar zI_m}\|=1.
$$
Clearly, 
$$
\|H_{\bar z^{d_0+1}\Psi+c\bar zI_m}\|_{\text e}<1.
$$
Let $G$ be a best approximation of $\bar z^{d_0+1}\Psi+c\bar zI_m$
by  $H^\be$ matrix functions. As we have explained in \S 2, 
$\bar z^{d_0+1}\Psi+c\bar zI_m-G$ admits a representation
\beq
\label{3.6}
\bar z^{d_0+1}\Psi+c\bar zI_m-G=W^*\left(\begin{array}{cc}u&0\\0&\U
\end{array}\right)V^*,
\end{equation}
where $V$ and $W^{\text t}$ are thematic matrices, $\|\U\|_{L^\be}\le1$,
$u$ is a unimodular function such that $T_u$ is Fredholm and $\ind T_u>0$.
Obviously, $\|H_{\bar z^{d_0+1}\Psi+c\bar zI_m-G}\|_{\text e}<1$. By Theorem 6.3 of
[PT], $\|H_\U\|_{\text e}<1$.

Let us show that $\ind T_u=1$. Suppose that $\ind T_u>1$. Then
$$
z(\bar z^{d_0+1}\Psi+c\bar zI_m-G)=W^*\left(\begin{array}{cc}zu&0\\0&z\U
\end{array}\right)V^*
$$
is still badly approximable (see \S 2). Hence,
$$
1=\|z(\bar z^{d_0+1}\Psi+c\bar zI_m-G)\|_{L^\be}=
\|H_{\bar z^{d_0}\Psi+cI_m-zG}\|=\|H_{\bar z^{d_0}\Psi}\|
\le\|\bar z^{d_0}\Psi\|_{L^\be}<1.
$$
We have got a contradiction. Multiplying both sides of (\ref{3.6}) by $z^{d_0+1}$,
we obtain
$$
\Psi+cz^{d_0}I_m-z^{d_0+1}G=
W^*\left(\begin{array}{cc}z^{d_0+1}u&0\\0&z^{d_0+1}\U
\end{array}\right)V^*.
$$
Put $u_0=z^{d_0+1}u$. Clearly, $\ind T_{u_0}=-d_0-1+\ind T_u=-d_0$.

Let us show that $\|H_{z^{d_0+1}\U}\|<1$. Consider the following factorization:
$$
\bar z^{d_0}\Psi+cI_m-zG=
W^*\left(\begin{array}{cc}zu&0\\0&z\U
\end{array}\right)V^*.
$$
Clearly, $\|H_{\bar z^{d_0}\Psi+cI_m-zG}\|=\|H_{\bar z^{d_0}\Psi}\|<1$
and $\|H_{z\U}\|_{\text e}=\|H_\U\|_{\text e}<1$.
By Lemma 4.4 of [AP], $\|H_{z\U}\|<1$. Hence,
$$
\|H_{z^{d_0+1}\U}\|\le\|H_{z\U}\|<1.
$$
We can apply now the induction 
hypotheses to $z^{d_0+1}\U$. Finally, By Lemma 1.5 of [PY1], 
there exists a function $F\in H^\be(\mm_{m,m})$ such $\Psi-F$ admits a desired
representation. $\bl$

\begin{lem}
\label{t3.6}
Let $U$ be an $n\times n$ matrix function of the form {\em(\ref{3.5})}, where
the $V_j$ and $W_j$ are as in Lemma \ref{3.5}, the $u_j$ are unimodular
functions such that the operators $T_{u_j}$ are Fredholm whose indices are 
arbitrary integers. If {\em$\|H_U\|_{\text e}<1$}, then {\em$\|H_{U^*}\|_{\text e}<1$}.
\end{lem}

{\bf Proof of Lemma \ref{t3.6}.} Since
$$
\|H_U\|_{\text e}=\|H_{z^lU}\|_{\text e},
\quad\|H_{U^*}\|_{\text e}=\|H_{z^lU^*}\|_{\text e}
$$
for any $l\in\Z$, we may assume without loss of generality that $\ind T_{u_j}>0$,
\linebreak$0\le j\le n-1$. In this case (\ref{3.5}) is a thematic factorization
of $U$. By Theorem 3.1 of [PY1], the Toeplitz operator $T_U$ has dense range
in $H^2(\C^n)$. (Theorem 3.1 is stated in [PY1] for $H^\be+C$ matrix functions
but it is easy to see that the proof given in [PY1] works for any functions that
admit thematic factorizations.) As in the proof of Lemma \ref{t1.1}, we can conclude 
that $\|H_U\|_{\text e}=\|H_{U^*}\|_{\text e}$ which proves the result. $\bl$

{\bf Proof of Theorem \ref{t3.3}.} If $r=n$ and $Q$ is a best approximation of
$\Phi$ by bounded analytic matrix functions, then $Q$ is a superoptimal approximation
of $\Phi$. It follows from (\ref{3.2}) that $\Phi-Q$ admits a thematic factorization
(see \S 1), and so $U=\Phi-Q$ is a unitary 
interpolant of $\Phi$. By Lemma \ref{t3.6}, $U$ satisfies (\ref{3.3}).

Suppose now that $r<n$ and $Q$ is a best approximation of $\Phi$ by bounded
analytic matrix functions. Then $\Phi-Q$ admits a factorization of the form
$$
\Phi-Q=W_0^*\cdots W^*_{r-1}
\left(\begin{array}{cccc}
u_0&\cdots&0&0\\
\vdots&\ddots&\vdots&\vdots\\
0&\cdots&u_{r-1}&0\\
0&\cdots&0&\Psi
\end{array}\right)V^*_{r-1}\cdots V_0^*,
$$
where
$$
W_j=\left(\begin{array}{cc}I_j&0\\0&\breve W_j\end{array}\right),\quad
V_j=\left(\begin{array}{cc}I_j&0\\0&\breve V_j\end{array}\right),\quad
1\le j\le r-1,
$$
$V_0$, $W_0^{\text t}$, $\breve V_j$, $\breve W^{\text t}_j$ are thematic 
matrix functions, $\|\Psi\|_{L^\be}\le1$, $\|H_\Psi\|<1$, $u_0,\cdots,u_{r-1}$ 
are unimodular functions such that the Toeplitz operators $T_{u_j}$ are Fredholm,
and
$$
\ind T_{u_0}\ge \ind T_{u_1}\ge\cdots\ge\ind T_{u_{r-1}}>0
$$
(see \S 2).

We can apply now Lemma \ref{t3.5} to $\Psi$ and find a matrix function 
$G\in H^\be(\mm_{n-r,n-r})$ such that
$$
\Psi-G=\check W_r^*\cdots\check W^*_{n-1}
\left(\begin{array}{ccc}
u_r&\cdots&0\\
\vdots&\ddots&\vdots\\
0&\cdots&u_{n-1}
\end{array}\right)\check V^*_{n-1}\cdots\check V_r^*,
$$
where
$$
\check W_j=\left(\begin{array}{cc}I_{j-r}&0\\0&\breve W_j\end{array}\right),\quad
\check V_j=\left(\begin{array}{cc}I_{j-r}&0\\0&\breve V_j\end{array}\right),\quad
r\le j\le n-1,
$$
$\breve V_j$, $\breve W^{\text t}_j$ are thematic matrix functions,
$u_r,\cdots,u_{n-1}$ are unimodular functions such that the operators $T_{u_j}$
are Fredholm and
$$
\ind T_{u_j}=-d_j,\quad r\le j\le n-1.
$$

Using the trivial part of Lemma 1.5 of [PY1], we can find inductively a matrix 
function $F\in H^\be(\mm_{n,n})$ such that $\Phi-F$ admits a factorization
$$
\Phi-F=W_0^*\cdots W^*_{n-1}
\left(\begin{array}{ccc}
u_0&\cdots&0\\
\vdots&\ddots&\vdots\\
0&\cdots&u_{n-1}
\end{array}\right)V^*_{n-1}\cdots V_0^*,
$$
where
$$
W_j=\left(\begin{array}{cc}I_j&0\\0&\breve W_j\end{array}\right),\quad
V_j=\left(\begin{array}{cc}I_j&0\\0&\breve V_j\end{array}\right),\quad
r\le j\le n-1.
$$

Clearly, $\Phi-F$ is a unitary-valued function. By Lemma \ref{t3.6},
it satisfies (\ref{3.3}). To prove that the Wiener--Hopf 
factorization indices of $\Phi-F$ are equal to
$$
-\ind T_{u_0},-\ind T_{u_1},\cdots,-\ind T_{u_{n-1}},
$$
it is sufficient to apply Theorem \ref{t3.2} to the matrix function
$\bar z^{d_{n-1}+1}(\Phi-F)$. $\bl$

{\bf Proof of Theorem \ref{t3.4}.} Suppose that all superoptimal singular values of 
$\Phi$ are equal to 1. Let $U$ be a unitary interpolant of $\Phi$. Clearly,
$\Phi-U$ is the superoptimal approximation of $\Phi$ which is unique because of 
conditions (\ref{3.1}) and (\ref{3.2}).

If $\Phi$ has superoptimal singular values less than 1, then by Theorem \ref{t3.3},
$\Phi$ has infinitely many unitary interpolants satisfying (\ref{3.3}). $\bl$

\

\setcounter{equation}{0}
\section{\bf Hereditary properties of unitary interpolants}

\

In this section we consider the following problem. Suppose that the initial 
matrix function $\Phi$ belongs to a certain function space $X$. The question
is whether one can obtain results similar to those of \S 3 for
unitary interpolants that belong to the same class $X$. We prove that this can
be done for two natural classes of function spaces. The first class consists of
so-called $\cR$-spaces (see [PK]). The second class of spaces consists of Banach 
algebras satisfying Axioms (A1)--(A4) below. In both cases we apply so-called
recovery theorems for unitary-valued functions obtained in [Pe3], [Pe4] and the
results of \S 3 of this paper.

It follows easily from Nehari's theorem (see e.g., [N], [PK]) that
for a function $\f\in L^2$ the Hankel operator $H_\f$
defined on the set of polynomials in $H^2$ extends to a bounded 
operator from $H^2$ to $H^2_-$ if and only if $\pp_-\f$ belongs to the space
$BMO$ of functions of bounded mean oscillation. Indeed, this
can easily be seen from the following description of $BMO$ due to C. Fefferman:
$$
BMO=\{f+\pp_+g:~f,g\in L^\be\}
$$
(see e.g., [G]). Similarly, it follows from Hartman's theorem
that $H_\f$ is compact if and only if $\pp_-\f$ belongs to the space
$VMO$ of functions of vanishing mean oscillation which can be seen from
the following description of VMO due to Sarason:
$$
VMO=\{f+\pp_+g:~f,g\in C(\T)\}
$$
(see e.g., [G]).

Therefore we can consider the problem of finding unitary interpolants
for functions $\Phi\in BMO(\mm_{n,n})$.

We are not going to give a precise definition of $\cR$-spaces and refer the reader
to [PK] for details. Roughly speaking $\cR$-spaces are function spaces $X\subset VMO$
such that one can determine whether a function $\f$ belongs to $X$ by the behavior of
the singular values of the Hankel operators $H_\f$ and $H_{\bar\f}$ 

Important examples of $\cR$-spaces are the Besov spaces $B_p^{1/p}$, $0<p<\be$,
(this follows from the fact that $H_\f$ belongs to the Schatten-von Neumann
class $\bS_p$ if and only if $\pp_-\f\in B_p^{1/p}$, see [Pe1], [Pe2], [S])
and the space $VMO$ of functions of vanishing mean oscillation (this follows
from the above compactness criterion).

For a function space $X$, $X\subset VMO$, and a matrix function $\Phi\in X(\mm_{n,n})$
we consider now the problem of finding unitary interpolants $U$ of $\Phi$
that belong to the same class $\Phi\in X(\mm_{n,n})$.

Note that if $\Phi\in VMO(\mm_{n,n})$, then
$H_\Phi$ is a compact operator from $H^2(\C^n)$ to $H^2_-(\C^n)$, and so
the results of \S 3 are applicable to $\Phi$.

In addition to the class of $\cR$-spaces we consider the class of function spaces
$X\subset C(\T)$ that contain the trigonometric polynomials 
and satisfy the following axioms:

(A1) {\it If $f\in X$, then $\bar f\in X$ and $\pp_+f\in X$;}

(A2) {\it$X$ is a Banach algebra with respect to pointwise multiplication;}

(A3) {\it for every $\f\in X$ the Hankel operator $H_\f$ is a compact operator
from $X_+$ to $X_-$;}

(A4) {\it if $f\in X$ and $f$ does not vanish on $\T$, then $1/f\in X$.}

Here we use the notation
$$
X_+=\{f\in X:~\hat f(j)=0,~j<0\},\quad X_-=\{f\in X:~\hat f(j)=0,~j\ge0\}.
$$

Note that in [PK] and [Pe4] a similar system of axioms was considered. Axioms 
(A1), (A2), and (A4) were stated there in the same way but (A3) was stated in the
following way:

(A$'$3) {\it the trigonometric polynomials are dense in $X$.} 

It is easy to see that (A1), (A2), and (A$'$3) imply (A3). 
In [PK] many examples of function spaces were given that satisfy (A1), (A2),
(A$'$3), and (A4). Let us mention here the Besov classes $B_{p,q}^s$,
$1\le p<\be$, $1\le q\le\be$, $s>1/p$, the space of functions with absolutely
convergent Fourier series, the spaces
$$
\{f:~f^{(n)}\in VMO\},\quad \{f:~\pp_+f^{(n)}\in C(\T),~\pp_-f^{(n)}\in C(\T)\},
\quad n\ge1.
$$

In [PK] and [Pe4] nonseparable function spaces were treated  in a different way.
On the other hand
in [Pe4] another sufficient condition was given that implies (A3) and 
many examples of nonseparable spaces of functions are found that satisfy (A1)--(A4).
We mention here the H\"{o}lder--Zygmund spaces $\L_\a$, $\a>0$, the spaces
$$
\{f:~f^{(n)}\in BMO\},\quad \{f:~\pp_+f^{(n)}\in H^\be,~\pp_-f^{(n)}\in H^\be\},
\quad n\ge1,
$$
the space
$$
\{f:~|\hat f(j)|\le\const(1+|j|)^{-\a}\},\quad\a>1.
$$

To find a unitary interpolant in $X$ for a given matrix function $\Phi$ in $X$,
we need the following fact:

{\it Let $X$ be either an $\cR$-space or a space of functions 
satisfying {\em(A1)--(A4)} and let $U$ be a unitary-valued matrix function in 
$X(\mm_{n,n})$ such that the Toeplitz operator $T_U$ on $H^2(\C^n)$ is Fredholm.
Then}
\beq
\label{4.1}
\pp_-U\in X\quad\Longrightarrow\quad U\in X.
\end{equation}

For $\cR$-spaces this was proved in [Pe3], for spaces satisfying (A1), (A2), (A$'$3),
and (A4) this was proved in [Pe4]. Another method was used in [Pe4] to treat 
nonseparable Banach spaces. Note that for scalar unimodular functions such results
were obtained earlier in [PK]. Note also that (\ref{4.1}) remains true if we replace 
the condition that $T_U$ is Fredholm with the condition that $T_U$ has dense 
range in $H^2(\C^n)$.

The method used in [PK] and [Pe4] relies on hereditary properties of maximizing 
vectors of Hankel operators, it was used in [AAK1] in the special case of the 
space of functions with absolutely converging Fourier series. 

It was shown in [Pe4] that to prove (\ref{4.1}) for spaces satisfying (A1)--(A4), 
it suffices to prove the following assertions. 

\begin{lem}
\label{t4.1}
Suppose that $X$ is a function space satisfying {\em(A1)--(A3)}.
Let $\Phi\in X(\mm_{m,n})$ and let $\f\in H^2(\C^n)$ be a maximizing vector
of the Hankel operator $H_\Phi:H^2(\C^n)\to H^2_-(\C^m)$. Then $\f\in X_+(\C^n)$.
\end{lem}

\begin{lem}
\label{t4.2}
Suppose that $X$, $\Phi$, and $\f$ satisfy the hypotheses of Lemma \ref{t4.1}.
If $\f(\t)=0$ for some $\t\in\T$, then $(1-\bar\t z)^{-1}\f\in X_+(\C^n)$ and 
$(1-\bar\t z)^{-1}\f$ is also a maximizing vector of $H_\Phi$.
\end{lem}

Lemmas \ref{t4.1} and \ref{t4.2} were proved in [PY3] for spaces satisfying
(A1), (A2), and (A$'$3). We adjust slightly the proof given there
for spaces satisfying (A1)--(A3). Note that in the scalar case such results
were obtained in [PK]. Similar results were obtained in [DG1] by different methods.

{\bf Proof of Lemma \ref{t4.1}.} Without loss of generality we may assume that
the norm of the Hankel operator $H_\Phi:H^2(\C^n)\to H^2_-(\C^m)$ is equal to 1.

Consider the self-adjoint operator 
$H^*_\Phi H_\Phi$ on $H^2(\C^n)$. It follows from (A1) and (A2) that
it maps $X_+(\C^n)$ into itself. Let $R$ be the restriction of $H^*_\Phi H_\Phi$
to $X_+(\C^n)$. By (A3), $R$ is a compact operator on $X_+(\C^n)$. 

Clearly, $X_+(\C^n)\subset H^2(\C^n)$. We can imbed naturally the space 
$H^2(\C^n)$ to the dual space
$X_+^*(\C^n)$ as follows. Let $g\in H^2(\C^n)$. We associate with it the linear 
functional ${\cal J}(g)$ on $X_+(\C^n)$ defined by:
$$
f\mapsto(f,g)=\int_\T (f(\z),g(\z))_{\C^n}d\m(\z).
$$
Note that ${\cal J}(\l_1g_1+\l_2g_2)=\ov{\l_1}{\cal J}(g_1)+\ov{\l_2}{\cal J}(g_2)$, 
$g_1,\,g_2\in H^2$, $\l_1,\,\l_2\in\C$.

The imbedding ${\cal J}$ allows us to identify $H^2(\C^n)$ with a linear subset of 
$X_+^*(\C^n)$. Since $H^*_\Phi H_\Phi$ is self-adjoint, it is easy to see that
$R^*g=H^*_\Phi H_\Phi g$ for $g\in H^2(\C^n)$. Hence, 
$$
\Ker(I-R)\subset\Ker(I-H^*_\Phi H_\Phi)\subset\Ker(I-R^*).
$$
Since $R$ is a compact operator, it follows from the Riesz--Schauder theorem
(see [Yo], Ch. X, \S 5) that $\dim\Ker(I-R)=\dim\Ker(I-R^*)$.
Therefore
\beq
\label{4.2}
\Ker(I-R)=\Ker(I-H^*_\Phi H_\Phi)=\Ker(I-R^*).
\end{equation}

Clearly, the subspace $\Ker(I-H^*_\Phi H_\Phi)$ is the space of maximizing
vectors of $H_\Phi$, and it follows from (\ref{4.2}) that each maximizing vector 
of $H_\Phi$ belongs to $X_+$. $\bl$

{\bf Proof of Lemma \ref{t4.2}.} Suppose now that $\f$ is a maximizing vector
of $H_\Phi$ and $\f(\t)=0$. Consider the following continuous 
linear functional $\o$ on $X_+(\C^n)$:
\beq
\label{4.3}
\o(f)=(\pp_+(\f^*f))(\t).
\end{equation}
Let us show that $R^*\o=\o$. We have
$$
(R^*\o)(f)=\o(Rf)=\o(H^*_\Phi H_\Phi f)=(\pp_+(\f^* H^*_\Phi H_\Phi f))(\t)
$$
by (\ref{4.3}). Therefore
\begin{eqnarray*}
(R^*\o)(f)&=&(\pp_+(\f^*\pp_+(\Phi^*H_\Phi f)))(\t)\\
&=&(\pp_+(\f^*\Phi^*H_\Phi f))(\t)-(\pp_+(\f^*\pp_-(\Phi^*H_\Phi f)))(\t).
\end{eqnarray*}
Since $\f^*\pp_-(\Phi^*H_\Phi f)\in H^2_-$, it follows that 
$\pp_+(\f^*\pp_-(\Phi^*H_\Phi f))=0$. Hence,
\begin{eqnarray*}
(R^*\o)(f)&=&(\pp_+(\f^*\Phi^*H_\Phi f))(\t)\\&=&
(\pp_+((\pp_-\Phi\f)^*H_\Phi f))(\t)+(\pp_+((\pp_+\Phi\f)^*H_\Phi f))(\t).
\end{eqnarray*}
Clearly, $\pp_+((\pp_+\Phi\f)^*H_\Phi f)=0$, and so
\begin{eqnarray*}
(R^*\o)(f)&=&(\pp_+((H_\Phi\f)^*H_\Phi f))(\t)=(\pp_+((H_\Phi\f)^*\pp_-(\Phi f)))(\t)\\
&=&(\pp_+((H_\Phi\f)^*\Phi f))(\t)-(\pp_+((H_\Phi\f)^*\pp_+(\Phi f)))(\t).
\end{eqnarray*}
Since $(H_\Phi\f)^*\pp_+(\Phi f)\in H^2$, we have
\beq
\label{4.4}
(\pp_+((H_\Phi\f)^*\pp_+(\Phi f)))(\t)=(H_\Phi\f)^*(\t)(\pp_+(\Phi f))(\t).
\end{equation}
It is well known (see e.g., Theorem 0.2 of [PY1]) that for a maximizing vector
$\f$ of $H_\Phi$ the following equality holds:
\beq
\label{4.5}
\|(H_\Phi\f)(\z)\|_{\C^m}=\|H_\Phi\|\cdot\|\f(\z)\|_{\C^n},\quad\z\in\T.
\end{equation}
By the hypotheses $\f(\t)=0$, and so by (\ref{4.5}), $(H_\Phi\f)(\t)=0$. Hence,
it follows from (\ref{4.4}) that
\begin{eqnarray*}
(R^*\o)(f)&=&(\pp_+((H_\Phi\f)^*\Phi f))(\t)=(\pp_+((\Phi^*H_\Phi\f)^*f))(\t)\\
&=&(\pp_+((\pp_+(\Phi^*H_\Phi\f))^*f))(\t)+(\pp_+((\pp_-(\Phi^*H_\Phi\f))^*f))(\t)\\
&=&(\pp_+((H^*_\Phi H_\Phi\f)^*f))(\t)+((\pp_-(\Phi^*H_\Phi\f)^*)f)(\t),
\end{eqnarray*}
since $(\pp_-(\Phi^*H_\Phi\f))^*\in H^2$. Since $\f$ is a maximizing vector of 
$H_\Phi$, it follows that $H^*_\Phi H_\Phi\f=\f$, and so
\begin{eqnarray*}
(R^*\o)(f)&=&(\pp_+(\f^*f))(\t)+((\pp_-(\Phi^*H_\Phi\f)^*)f)(\t)\\
&=&\o(f)+(\Phi^*H_\Phi\f)^*(\t)f(\t)-((\pp_+(\Phi^*H_\Phi\f)^*)f)(\t)\\
&=&\o(f)+\Phi^*(\t)((H_\Phi\f)(\t))^*f(\t)-(\f(\t))^*f(\t)=\o(f)
\end{eqnarray*}
which proves (\ref{4.3}).

It follows from (\ref{4.3}) that $\o\in\Ker(I-R^*)$, and by Lemma \ref{t4.1},
$\o\in\Ker(I-R)$. Hence, there must be a function $\psi\in X_+(\C^n)$ such that
$\psi$ is a maximizing vector of $H_\Phi$ and
$$
\o(f)=\int_\T(f(\z),\psi(\z))_{\C^n}d\m(\z),\quad f\in X_+(\C^n).
$$
To prove that $\psi=(1-\bar\t z)^{-1}\f$, we show that these two functions have the
same Taylor coefficients. Let $k\in\Z_+$, $\g\in\C^n$. We have
$$
(\hat\psi(k),\g)_{\C^n}=\ov{\o(z^k\g)}=\ov{(\pp_+z^k\f^*\g)(\t)}
=\ov{\sum_{j=0}^k\bar{\t}^{k-j}(\hat\f(j))^*\g},
$$
and so
$$
\hat\psi(k)=\sum_{j=0}^k\bar{\t}^{k-j}\hat\f(j).
$$
On the other hand
$$
(1-\bar\t z)^{-1}\f=\left(\sum_{k\ge0}\bar\t^kz^k\right)\f=
\sum_{k\ge0}z^k\left(\sum_{j=0}^k\bar\t^{k-j}\hat\f(j)\right)
$$
which completes the proof. $\bl$

Now we are ready to state the main result of this section. Recall that if $U$ is
a unitary-valued matrix function in $VMO$, then $H_U$ and $H_{U^*}$ are compact,
and so the Toeplitz operator $T_U$ is Fredholm (see \S 2).

\begin{thm}
\label{t4.3}
Suppose that $X$ is either an $\cR$-space or a space of functions satisfying
{\em(A1)--(A4)}. Let $\Phi$ be a matrix function in $X(\mm_{n,n})$ such that
$\|H_\Phi\|\le1$. Let $r$ be the number of superoptimal singular values of $\Phi$
equal to 1. Then for any integers $\{d_j\}_{r\le j\le n-1}$ such that
$$
0\le d_r\le d_{r+1}\le\cdots\le d_{n-1}
$$
there exists a unitary interpolant $U\in X(\mm_{n,n})$ whose nonnegative
Wiener--Hopf factorization indices are $d_r,d_{r+1},\cdots,d_{n-1}$.
\end{thm}

Note that any unitary interpolant $U\in X(\mm_{n,n})$ must satisfy (\ref{3.3}), and 
so by Theorem \ref{t3.2}, $U$ must have precisely $r$ negative Wiener--Hopf
factorization indices that are uniquely determined by $\Phi$.

{\bf Proof of Theorem \ref{t4.3}.} By Theorems \ref{t3.2} and \ref{t3.3},
it is sufficient to prove that under the hypotheses
of the theorem any unitary interpolant $U$ of $\Phi$ that satisfies (\ref{3.3}) 
must belong to $X$. 

As we have already observed in the introduction, if $U$ is a unitary interpolant
of $\Phi$ that satisfies (\ref{3.3}), then
$T_U$ is Fredholm.  Hence, the desired result is just
implication (\ref{4.1}) which is true both for the $\cR$-spaces $X$ (see [Pe3]) and
for the spaces $X$ satisfying (A1)--(A4) (this follows from Lemmas \ref{t4.1}, 
\ref{t4.2} and the results of [Pe4]). $\bl$

Note that Theorem \ref{t4.3} implies the main results (in particular, Theorem 1.1)
stated in [DG1] .

The following special case of Theorem \ref{t4.3} is particularly important.

\begin{thm}
\label{t4.4}
Let $\Phi$ be a matrix function in $(H^\be+C)(\mm_{n,n})$ such that \linebreak
$\|H_\Phi\|\le1$. Let $r$ be the number of superoptimal singular values of $\Phi$
equal to 1. Then for any integers $\{d_j\}_{r\le j\le n-1}$ such that
$$
0\le d_r\le d_{r+1}\le\cdots\le d_{n-1}
$$
there exists a unitary interpolant $U\in QC(\mm_{n,n})$ whose nonnegative
Wiener--Hopf factorization indices are $d_r,d_{r+1},\cdots,d_{n-1}$.
\end{thm}

\Pf Let $X=VMO$. As we have already observed, $X$ is an $\cR$-space.
The condition $\Phi\in H^\be+C$ implies $\pp_-\Phi\in X$. Let $U$ be a unitary 
interpolant of $\Phi$ that satisfies the conclusion of Theorem \ref{t3.3}.
By Theorem \ref{t4.3}, $U\in VMO$. The result follows now from the well-known
identity
$$
QC=VMO\cap L^\be
$$
(see e.g., [G]). $\bl$

\

\noindent
Department of Mathematics
\newline
Kansas State University
\newline
Manhattan, Kansas 66506
\newline
USA


\begin{thebibliography}{99}

\bibitem[AAK1]{AAK1} {\sc V.M. Adamyan, D.Z. Arov, and M.G. Krein},
On infinite Hankel matrices and generalized problems of 
Carath\'{e}odory-Fej\'{e}r and F. Riesz, 
{\em Funktsional. Anal. i Prilozhen.}
{\bf 2:1} (1968), 1-19 (In Russian).

\bibitem[AAK2]{AAK2} {\sc V.M. Adamyan, D.Z. Arov, and M.G. Krein},
On infinite Hankel matrices and generalized problems of
Carath\'{e}odory-Fej\'{e}r and I. Schur,
{\em Funktsional. Anal. i Prilozhen.}
{\bf 2:2} (1968), 1-17 (In Russian).

\bibitem[AP]{AP} {\sc R.B. Alexeev and V.V. Peller}, Invariance properties of thematic
factorizations of matrix functions, to appear

\bibitem[B]{B} {\sc J.A. Ball}, Invariant subspace representations, unitary 
interpolants and factorization indices, in ``Topics in Operator Theory: Advances
and Applications'' {\bf12} (1984), 11-38.

\bibitem[D]{D} {\sc R.G. Douglas}, ``Banach algebra techniques in operator theory'',
Academic Press, New York--London 1972.

\bibitem[DG1]{DG1} {\sc H. Dym and I. Gohberg}, Unitary interpolants, factorization
indices and infinite Hankel block matrices, 
{\em J. Funct. Anal.} {\bf54} (1983), 229-289.

\bibitem[DG2]{DG2} {\sc H. Dym and I. Gohberg}, On unitary interpolants and Fredholm
infinite block Toeplitz matrices, {\em Int. Eq. Op. Theory} {\bf6} (1983), 863-878.

\bibitem[G]{G} {\sc J. Garnett}, ``Bounded analytic functions'',
Academic Press, NY--London--Toronto--Sydney--San Franciso, 1981.

\bibitem[Kh]{Kh} {\sc S. Khavinson}, On some extremal problems of the theory 
of analytic functions, {\em Uchen. Zapiski Mosk. Universiteta, Matem.} {\bf 144}:4
(1951), 133-143. English Translation: Amer. Math. Soc. Translations (2) {\bf 32}
(1963), 139-154.

\bibitem[LS]{LS} {\sc G.S. Litvinchuk} and {\sc I.M. Spitkovski}, 
``Factorization of measurable 
matrix functions'', Oper. Theory: Advances and 
Appl., {\bf25}. Birkh\"{a}user Verlag, Basel-Boston, MA, 1987.

\bibitem[N]{N} {\sc N.K. Nikol'skii}, ``Treatise on the shift operator. Spectral
function theory,'' Springer-Verlag, Berlin--Heidelberg--New York--Tokyo, 1986.

\bibitem[Pa]{Pa} {\sc L. B. Page}, 
Bounded and compact vectorial Hankel operators,
{\em Trans. Amer. Math. Soc.}, {\bf 150} (1970) 529-539.

\bibitem[Pe1]{Pe1} {\sc V. V. Peller}, 
Hankel operators of class $S_p$, and applications (rational
approximation, Gaussian processes, the majorant problem for operators),
{\em Mat. Sbornik,} {\bf 113} (1980), 538-581. English Translation;
{\em Math USSR Sbornik,} {\bf 41} (1982), 443-479.

\bibitem[Pe2]{Pe2}{\sc V.V. Peller}, A description of Hankel operators of class
${\frak S}_p$ for $p>0$, investigation of the rate of rational approximation
and other applications, {\em Mat. Sb.} {\bf122} (1983), 481-510;
English Translation: {\em Math. USSR-Sb.} {\bf50} (1985), 465-494.

\bibitem[Pe3]{Pe3}{\sc V.V. Peller}, Hankel operators and multivariate stationary
processes, {\em Proc. Symp. Pure Math.} {\bf51} (1990), 357-371.

\bibitem[Pe4]{Pe4}{\sc V.V. Peller}, Factorization and approximation 
problems for matrix functions, {\em Amer. J. Math.} {\bf11} (1998), 751-770.

\bibitem[PK]{PK} {\sc V.V. Peller} and {\sc S.V. Khruschev}, Hankel operators, 
best approximation and stationary Gaussian processes, 
{\em Uspekhi Mat. Nauk} {\bf 37:1} (1982), 53-124.
English Translation: {\em Russian Math. Surveys} {\bf 37} (1982),
53--124.

\bibitem[PT]{PT} {\sc V.V. Peller and S.R. Treil}, 
Approximation by analytic matrix functions.
The four block problem, {\em J. Funct. Anal.} {\bf148} (1997), 191-228.

\bibitem[PY1]{PY1} {\sc V.V. Peller and N.J. Young}, Superoptimal analytic 
approximations of 
matrix functions, {\em J. Funct. Anal.} {\bf 120} (1994), 300-343. 

\bibitem[PY2]{PY2} {\sc V.V. Peller and N.J. Young},
Superoptimal singular values and indices of
matrix functions, {\em Int. Eq. Op. Theory} {\bf20} (1994), 35-363.

\bibitem[PY3]{PY3} {\sc V.V. Peller and N.J. Young},
Continuity properties of best analytic approximation,
{\em J. Reine Angew. Math.} {\bf483} (1997), 1-22.

\bibitem[Po]{Po} {\sc S. J. Poreda}, 
A characterization of badly approximable functions,
{\em Trans. Amer. Math. Soc.}, {\bf 169} (1972), 249-256.

\bibitem[Sa]{Sa} D. Sarason, ``Function theory on the unit circle,''
Notes for lectures at Virginia Polytechnic Inst. and State Univ., 1978.

\bibitem[Se]{Se} {\sc S. Semmes}, Trace ideal criteria for Hankel operators and 
applications to Besov spaces, {\em Integral Equations and Operator Theory}, 
{\bf 7} (1984), 241-281.

\bibitem[Si]{Si} {\sc I.B. Simonenko}, Some general problems of the theory of
the Riemann boundary-value problem, {\em Izv. Akad. Nauk SSSR, Ser. Mat.}
{\bf32} (1968), 1138-1146 (In Russian).

\bibitem[T]{T} {\sc S.R. Treil}, On superoptimal approximation by analytic and 
meromorphic matrix-valued functions, 
{\em J. Functional Analysis} {\bf 131} (1995), 386-414.

\bibitem[Yo]{Yo} {\sc K. Yosida}, ``Functional Analysis" (Sixth edition), 
Springer Verlag, Berlin, 1980.



\end{thebibliography}
\end{document}